\documentclass[twoside,reqno]{amsart} 

\usepackage{latexsym,rotate,eucal,cite,url}
\usepackage{amsmath,amsthm,amssymb,amsxtra} 
\usepackage{enumitem}
\usepackage{color,graphics,graphicx}

\numberwithin{subsection}{section}

\newcommand{\puad}{\hspace*{-3.5mm}}

\newcommand{\fns}{\footnotesize}

\newcommand{\sof}[1]{\lceil{#1}\rceil}
\newcommand{\pav}[1]{\lfloor{#1}\rfloor}
\newcommand{\soff}[1]{\left\lceil{#1}\right\rceil}
\newcommand{\pavv}[1]{\left\lfloor{#1}\right\rfloor}

\newcommand{\res}[1]{\underset{#1}{\md{Res}}}

\newcommand{\mb}[1]{\mathbb{#1}}
\newcommand{\mc}[1]{\mathcal{#1}}
\newcommand{\md}{\mathrm}

\newcommand{\sbs}[1]{_{\substack{#1}}}

\newcommand{\binm}{\binom}

\newcommand{\be}{\begin{equation}}
\newcommand{\ee}{\end{equation}}
\newcommand{\ba}{\begin{array}}
\newcommand{\ea}{\end{array}}
\newcommand{\bmn}{\begin{eqnarray}}
\newcommand{\emn}{\end{eqnarray}}
\newcommand{\bnm}{\begin{eqnarray*}}
\newcommand{\enm}{\end{eqnarray*}}
\newcommand{\bln}{\begin{subequations}}
\newcommand{\eln}{\end{subequations}}

\newcommand{\pq}[1]{\begin{equation}#1\end{equation}}
\newcommand{\pnq}[1]{\begin{align*}#1
            \end{align*}}    

\newcommand{\centro}[1]
           {\begin{center}#1\end{center}}

\newcommand{\lam}{\lambda}

\newcommand{\veps}{\varepsilon}
\newcommand{\del}{\delta}

\newtheorem{thm}{Theorem}
\newtheorem{lemm}[thm]{Lemma}

\newtheorem{entry}{Entry}

\newcommand{\referxy}[4]{\bibitem{kn:#1}{#2,}~\emph{#3,}~{#4.}}	
\newcommand{\cito}[1]{\cite{kn:#1}}	
\newcommand{\citu}[2]{\cite[#2]{kn:#1}}

\begin{document}{\fbox{\fns\today}\hfill}
\title{Integrals of Hyperbolic Tangent Function}
\author{Jing Li and Wenchang Chu}
\dedicatory{\emph{
School of Mathematics and Statistics,
Zhoukou Normal University, Henan, China\\
Department of Mathematics and Physics,
University of Salento, Lecce 73100, Italy}}
\thanks{Corresponding author (W.~Chu): chu.wenchang@unisalento.it}

\begin{abstract}
By means of the contour integration method, we evaluate, in closed form,
a class of definite integrals involving hyperbolic tangent function.
\end{abstract}

\maketitle\thispagestyle{empty}\vspace*{-5mm} 
\markboth{Jing Li and Wenchang Chu}
         {Integrals of Hyperbolic Tangent Function}


There exist numerous interesting and challenging integrals in the
mathematical literature (see for example \cite{kn:boros,kn:valean}).
However, the following improper integrals don't seem to have been
examined previously
\pq{J(m,n):=\label{int=mn}
\int^{\infty}_{0}\frac{\tanh^m(z)}{z^n}dz,}
where $m,n\in\mb{N}$ subject to $m\ge n\ge 2$ and $m\equiv_{2} n$.
Here and forth $m\equiv_{2} n$ stands for that $m$ is congruent
to $n$ modulo 2. The aim of this short paper is to establish
explicit formula for $J(m,n)$ when $m$ and $n$ are considered
as two integer parameters. The final result states that the
integral value $J(m,n)$ is expressed as a finitely linear
combination of $\frac{\zeta(m+n-2k+1)}{\pi^{m+n-2k}}$ for
$1\le k\le \soff{\frac{m}2}$, where $\sof{x}$ denotes the
smallest integer not less than $x$.
This will be fulfilled by employing
the contour integration method.

\subsection{The contour integral} \
Rewriting the integrand in terms of the exponential function
\pq{\label{T-mn}T_{m,n}(z):=\frac{\tanh^m(z)}{z^n}=\frac{1}{z^n}
\bigg(\frac{e^z-e^{-z}}{e^z+e^{-z}}\bigg)^m
=\frac{1}{z^{n}}\bigg(\frac{e^{2z}-1}{e^{2z}+1}\bigg)^{m},}
we can express, by symmetry, the integral as
\[2J(m,n)=2\int^{\infty}_{0}\frac{1}{z^{n}}
\bigg(\frac{e^{2z}-1}{e^{2z}+1}\bigg)^{m}dz
=\int^{\infty}_{-\infty}\frac{1}{z^{n}}
\bigg(\frac{e^{2z}-1}{e^{2z}+1}\bigg)^{m}dz.\]
Observe that $T_{m,n}(z)$ is a meromorphic function
whose poles of $T_{m,n}(z)$ are determined by
\[1+e^{2z}=0.\]
We find that all the poles of $T_{m,n}(z)$ are given by
\[z_{k}=\bigg(k-\frac{1}{2}\bigg)\pi i,\quad\text{where}\qquad{k\in\mb{Z}}.\]

In order to evaluate the integral $J(m,n)$, consider the anti-clockwise
contour $\mc{L}$, consisting of two segments $[-\pi M,-\veps]$ and $[\veps,\pi M]$
along the real axis plus two semi-circles $\mc{C}(\veps)$ and $\mc{C}(\pi M)$ over
the upper half-plane centered at the origin with radii $\veps$ and $\pi M$
(where $M\in\mb{N}$), respectively.
Hence, all the poles of $T_{m,n}(z)$ inside the contour $\mc{L}$
are of order $m$ and given explicitly by
\[z_{k}=\bigg(k-\frac{1}{2}\bigg)\pi i,
\quad\text{where}\quad{k\in\mb{N}~\text{ with }~1\le k\le M}.\]

Then according to the residue theorem (see \cito{chu07q}
and \citu{titch}{\S3.3}), we have the following equality
\pnq{\int_{\mc{L}}T_{m,n}(z)dz=2\pi i\sum_{k=1}^{M}\res{z=z_{k}}T_{m,n}(z)
&=\int_{\mc{C}(\veps)}T_{m,n}(z)dz+\int_{\mc{C}(\pi M)}T_{m,n}(z)dz\\
&+\int^{-\veps}_{-\pi M}T_{m,n}(z)dz+\int^{\pi M}_{\veps}T_{m,n}(z)dz.}

As we shall show in the next two sections that
\[\lim_{\veps\to0}\int_{\mc{C}(\veps)}T_{m,n}(z)dz
=\lim_{M\to\infty}\int_{\mc{C}(\pi M)}T_{m,n}(z)dz=0,\]
the limiting case of the above equation as $\veps\to0$ and $M\to\infty$
will reduce to the following identity.
\begin{lemm}\label{int=res}
For two integers $m,n\in\mb{N}$ subject to conditions $m\ge n\ge 2$
and $m\equiv_2 n$, the following integral-sum formula holds:
\[J(m,n)=\int^{\infty}_{0}\frac{\tanh^m(z)}{z^n}dz
=\pi i\sum_{k=1}^{\infty}\res{z=z_{k}}T_{m,n}(z).\]
\end{lemm}

\subsection{The integral along $\mc{C}(\veps)$} \
Keeping in mind $m\ge n$, we have
\[\lim_{z\to0}T_{m,n}(z)
=\lim_{z\to0}\frac{1}{z^{n}}\bigg(\frac{e^{2z}-1}{e^{2z}+1}\bigg)^{m}
=\begin{cases}
1,&m=n,\\
0.&m>n.
\end{cases}\]
Therefore, for any $\veps$ with $|\veps|$ being sufficiently small,
the function $T_{m,n}(z)$ is bounded on $C_\veps$, which implies that
\[\bigg|\int_{\mc{C}(\veps)}T_{m,n}(z)dz\bigg|
=\mc{O}(1)\times\Big|\int_{\mc{C}(\veps)}dz\Big|
=\mc{O}(\veps\pi),\]
and consequently, the following limiting relation
\[\lim_{\veps\to0}\int_{\mc{C}(\veps)}T_{m,n}(z)dz=0.\]

\subsection{The integral along $\mc{C}(\pi M)$} \
Writing $z=\pi Me^{i\theta}$ on the semi-circle $\mc{C}(\pi M)$,
we can explicitly express the modulus
\[\big|1+e^{2z}\big|^2=\Big|1+e^{2\pi M(\cos\theta+i\sin\theta)}\Big|^2
=1+e^{4\pi M\cos\theta}+2e^{2\pi M\cos\theta}\cos(2\pi M\sin\theta).\]
Observe that the squared modulus results in 4 for $\theta=\frac{\pi}2$,
and is equal to or greater than $(1-e^{2\pi M\cos\theta})^2$ when
$\theta\ne\frac{\pi}2$. Then we can choose a small $\del>0$ such that
\[\big|1+e^{2z}\big|^2\ge\begin{cases}
3,&\frac{\pi}2-\del\le\theta\le \frac{\pi}2+\del;\\
(1-e^{-2\pi M\sin\del})^2,&0\le\theta<\frac{\pi}2-\del
\text{ and }\frac{\pi}2+\del<\theta\le \pi.
\end{cases}\]
Therefore $(1+e^{2z})^{-1}$ is bounded on $\mc{C}(\pi M)$,
which implies further that the function below is bounded too on $\mc{C}(\pi M)$
\[\bigg(\frac{e^{2z}-1}{e^{2z}+1}\bigg)^{m}=\bigg(1-\frac{2}{e^{2z}+1}\bigg)^{m}
=\mc{O}(1).\]
Under the change in variable $z=\pi Me^{i\theta}$,
we estimate the integral
\[\Big|\int_{\mc{C}(\pi M)}T_{m,n}(z)dz\Big|
=\mc{O}(1)\times\bigg|\int_0^{\pi}
\frac{ie^{(1-n)i\theta}}{(\pi M)^{n-1}}d\theta\bigg|
=\mc{O}\Big(\frac{\pi^{2-n}}{M^{n-1}}\Big)\]
Hence, we have shown that for $m\ge n\ge 2$:
\[\lim_{M\to\infty}\int_{\mc{C}(\pi M)}T_{m,n}(z)dz=0.\]

\subsection{Computing residues} \
Finally, we have to determine the residues of $T_{m,n}(z)$ at $z_{k}$ explicitly.
When $m\ge n\ge2$ are small integers, it is not difficult to do this. However,
it becomes quite a tough task if $m$ and $n$ are considered as integer parameters.

Instead of calculating higher derivatives, we shall determine the residues
by extracting coefficients from formal power series (cf.~\citu{comtet}{\S3.2}
and \citu{wilf}{\S2.1}).
For the sake of brevity, denote by $[x^n]\phi(x)$ the coefficient
of $x^n$ in the formal power series $\phi(x)$. Keeping in mind that
$e^{2z_k}=(-1)^{2k-1}=-1$, we can write
\pnq{\res{z=z_{k}}T_{m,n}(z)
&=\big[(z-z_{k})^{-1}\big]\frac{1}{z^{n}}
\bigg(\frac{e^{2z}-1}{e^{2z}+1}\bigg)^{m}\\
&=\big[(z-z_{k})^{-1}\big]\frac{1}{z^{n}}
\bigg(\frac{e^{2z}+e^{2z_{k}}}{e^{2z}-e^{2z_{k}}}\bigg)^{m}\\
&=\big[z-z_{k})^{-1}\big]\frac{1}{z^n}
\bigg(\frac{e^{z-z_{k}}+e^{z_{k}-z}}{e^{z-z_{k}}-e^{z_{k}-z}}\bigg)^{m}.}
By making the change of variables $y=z-z_k$, we have further
\pnq{\res{z=z_{k}}T_{m,n}(z)
&=\big[y^{-1}\big]\frac1{(y+z_{k})^{n}}
\bigg(\frac{e^{y}+e^{-y}}{e^{y}-e^{-y}}\bigg)^{m}\\
&=\big[y^{-1}\big]\frac{(e^{y}+e^{-y})^{m}}{(2y)^{m}(y+z_{k})^{n}}
\bigg(\frac{e^{y}-e^{-y}}{2y}\bigg)^{-m}.}
Then the above residue can be expressed as
\pq{\label{res=cc}\res{z=z_{k}}T_{m,n}(z)
=\big[y^{m-1}\big]\frac{\md{U}(y)\md{V}(y)}{(y+z_{k})^n}
=\sum\sbs{1\le j\le m\\j\:\equiv_2\:m}
\binm{-n}{j-1}\big[y^{m-j}\big]\frac{\md{U}(y)\md{V}(y)}{z_{k}^{n+j-1}},}
where $\md{U}(y)$ and $\md{V}(y)$ are two even functions and given by
\[\md{U}(y)=\frac{(e^{y}+e^{-y})^{m}}{2^m}
\quad\text{and}\quad
\md{V}(y)=\bigg\{1-\bigg(1-\frac{e^{y}-e^{-y}}{2y}\bigg)\bigg\}^{-m}.\]

In order to extract the coefficient of $y^{m-j}$ in \eqref{res=cc},
it is sufficient to expand $\md{U}(y)$ and $\md{V}(y)$ on the right
to Maclaurin polynomials up to order $m$.

First by means of the binomial theorem, $\md{U}(y)$ can be written as
\pq{\label{exp-u}
\md{U}(y)=\Big(\frac{e^{y}+e^{-y}}{2}\Big)^{m}
=\sum_{\lam=0}^{m}\binm{m}{\lam}\frac{e^{(m-2\lam)y}}{2^{m}}.}
Then by observing the initial terms of the Maclaurin series
\[\frac{e^{y}-e^{-y}}{2y}=1+\frac{y^2}{6}+\frac{y^4}{120}+\mc{O}(y^6),\]
we can expand $\md{V}(y)$ in succession as follows:
\pnq{\md{V}(y)
&=\mc{O}(y^m)+\sum_{\ell=0}^{\pav{\frac{m}{2}}}(-1)^{\ell}\binm{-m}{\ell}
\bigg(1-\frac{e^{y}-e^{-y}}{2y}\bigg)^{\ell}\\
&=\mc{O}(y^m)+\sum_{\ell=0}^{\pav{\frac{m}{2}}}(-1)^{\ell}\binm{-m}{\ell}
\sum_{\mu=0}^{\ell}(-1)^{\mu}\binm{\ell}{\mu}
\bigg(\frac{e^{y}-e^{-y}}{2y}\bigg)^{\mu}.}
Let $\md{S}$ stand for the above double sum.
According to the binomial relations
\pnq{\binm{-m}{\ell}\binm{\ell}{\mu}
&=\binm{-m}{\mu}\binm{-m-\mu}{\ell-\mu},\\
\sum_{\ell=\mu}^{\pav{\frac{m}{2}}}
(-1)^{\ell-\mu}\binm{-m-\mu}{\ell-\mu}
&=(-1)^{\pav{\frac{m}{2}}-\mu}
\binm{-m-\mu-1}{\pavv{\frac{m}{2}}-\mu};}
we can manipulate $\md{S}$ in the following manner:
\pnq{\md{S}&=\sum_{0\le\mu\le \ell\le\pav{\frac{m}{2}}}
(-1)^{\ell-\mu}\binm{-m}{\ell}\binm{\ell}{\mu}
\bigg(\frac{e^{y}-e^{-y}}{2y}\bigg)^{\mu}\\
&=\sum_{0\le\mu\le \ell\le\pav{\frac{m}{2}}}
(-1)^{\ell-\mu}\binm{-m}{\mu}\binm{-m-\mu}{\ell-\mu}
\bigg(\frac{e^{y}-e^{-y}}{2y}\bigg)^{\mu}\\
&=\sum_{\mu=0}^{\pav{\frac{m}2}}
\binm{-m}{\mu}\bigg(\frac{e^{y}-e^{-y}}{2y}\bigg)^{\mu}
\sum_{\ell=\mu}^{\pav{\frac{m}2}}
(-1)^{\ell-\mu}\binm{-m-\mu}{\ell-\mu}\\
&=\sum_{\mu=0}^{\pav{\frac{m}2}}
(-1)^{\pav{\frac{m}{2}}-\mu}\binm{-m}{\mu}
\binm{-m-\mu-1}{\pavv{\frac{m}{2}}-\mu}
\bigg(\frac{e^{y}-e^{-y}}{2y}\bigg)^{\mu}.}
Applying further the binomial theorem to the rightmost function
and then making substitutions, we derive the following double
sum expansion
\pq{\label{exp-v}
\md{V}(y)=\mc{O}(y^m)+\puad
\sum_{0\le\nu\le\mu\le\pav{\frac{m}{2}}}
(-1)^{\pav{\frac{m}{2}}-\mu-\nu}
\binm{-m-\mu-1}{\pavv{\frac{m}{2}}-\mu}
\binm{-m}{\mu}\binm{\mu}{\nu}
\frac{e^{(\mu-2\nu)y}}{(2y)^{\mu}}.}

Now substituting \eqref{exp-u} and \eqref{exp-v} into \eqref{res=cc},
we can  further reformulate the residue as
\pnq{\res{z=z_{k}}T_{m,n}(z)
=&\sum\sbs{1\le j\le m\\j\:\equiv_2\:m}\sum_{\lam=0}^{m}
\binm{-n}{j-1}\binm{m}{\lam}
\sum_{0\le\nu\le\mu\le\pav{\frac{m}{2}}}
\binm{-m-\mu-1}{\pavv{\frac{m}{2}}-\mu}\\
&\times
\frac{(-1)^{\pav{\frac{m}{2}}-\mu-\nu}}{2^{m+\mu}z_{k}^{n+j-1}}
\binm{-m}{\mu}\binm{\mu}{\nu}
\big[y^{m-j}\big]\frac{e^{(m+\mu-2\nu-2\lam)y}}{y^{\mu}}\\
=&\sum\sbs{1\le j\le m\\j\:\equiv_2\:m}\sum_{\lam=0}^{m}
\binm{-n}{j-1}\binm{m}{\lam}
\sum_{0\le\nu\le\mu\le\pav{\frac{m}{2}}}
\binm{-m-\mu-1}{\pavv{\frac{m}{2}}-\mu}\\
&\times
\frac{(-1)^{\pav{\frac{m}{2}}-\mu-\nu}}{z_{k}^{n+j-1}}
\binm{-m}{\mu}\binm{\mu}{\nu}
\frac{(m+\mu-2\nu-2\lam)^{m+\mu-j}}{2^{m+\mu}(m+\mu-j)!}.}
Taking into account the binomial equality
\pnq{\binm{-m}{\mu}\binm{-m-\mu-1}{\pavv{\frac{m}{2}}-\mu}
&=(-1)^{\pavv{\frac{m}2}}\frac{(m)_{\mu}}{\mu!}
\frac{(1+m+\mu)_{\pavv{\frac{m}2}-\mu}}{(\pavv{\frac{m}2}-\mu)!}\\
&=(-1)^{\pavv{\frac{m}2}}\binm{\pavv{\frac{m}2}}{\mu}
\binm{m+\pavv{\frac{m}2}}{m}\frac{m}{m+\mu},}
we find finally the following quadruplicate sum expression
\pnq{\res{z=z_{k}}T_{m,n}(z)
=&\binm{m+\pavv{\frac{m}2}}{m}
\sum\sbs{1\le j\le m\\j\:\equiv_2\:m}
\binm{-n}{j-1}
\sum_{\lam=0}^{m}\sum_{\mu=0}^{\pav{\frac{m}2}}\sum_{\nu=0}^{\mu}
\frac{(-1)^{\mu+\nu}}{z_{k}^{n+j-1}}\\
&\times
\frac{m}{m+\mu}
\binm{m}{\lam}\binm{\pavv{\frac{m}2}}{\mu}\binm{\mu}{\nu}
\frac{(m-2\lam+\mu-2\nu)^{m+\mu-j}}{2^{m+\mu}(m+\mu-j)!}.}

\subsection{Conclusive theorem} \
To evaluate the sum of residues in Lemma~\ref{int=res}, it is enough to do that for
\[\sum_{k=1}^{\infty}\frac1{z_{k}^{n+j-1}}
=\bigg(\frac{2}{\pi i}\bigg)^{n+j-1}
\sum_{k=1}^{\infty}\frac1{(2k-1)^{n+j-1}}
=\frac{2^{n+j-1}-1}{(\pi i)^{n+j-1}}
\zeta(n+j-1).\]
Consequently, we have proved the following general integral identity.
\begin{thm}
Let $m$ and $n$ be the two natural numbers of the same parity with $m\ge n\ge2$,
the following integral formula holds:
\pnq{J(m,n)&=\int^{\infty}_{0}\frac{\tanh^m(z)}{z^n}dz
=\binm{m+\pavv{\frac{m}2}}{m}
\sum\sbs{1\le j\le m\\j\:\equiv_2\:m}
\binm{-n}{j-1}\frac{2^{n+j-1}-1}{(\pi i)^{n+j-2}}
\zeta(n+j-1)\\
&\times
\sum_{\lam=0}^{m}\sum_{\mu=0}^{\pav{\frac{m}2}}
\sum_{\nu=0}^{\mu}(-1)^{\mu+\nu}\frac{m}{m+\mu}
\binm{m}{\lam}\binm{\pavv{\frac{m}2}}{\mu}\binm{\mu}{\nu}
\frac{(m-2\lam+\mu-2\nu)^{m+\mu-j}}{2^{m+\mu}(m+\mu-j)!}.}
\end{thm}
From this theorem, we assert that for each pair of $m,n\in\mb{N}$
subject to $m\equiv_2n$ and $m\ge n\ge2$, the integral value
$J(m,n)$ results always in a finitely linear combination
of $\frac{\zeta(m+n-2k+1)}{\pi^{m+n-2k}}$ with
$1\le k\le \soff{\frac{m}2}$. The initial values
of these integrals are recorded as follows:

\pnq{J(2,2)&=\frac{14 \zeta (3)}{\pi ^2},\\
J(3,3)&=-\frac{7 \zeta (3)}{\pi ^2}+\frac{186 \zeta (5)}{\pi ^4},\\
J(4,2)&=\frac{56 \zeta(3)}{3 \pi ^2}-\frac{124 \zeta (5)}{\pi ^4},\\
J(4,4)&=-\frac{496 \zeta (5)}{3 \pi ^4}+\frac{2540 \zeta (7)}{\pi ^6},\\
J(5,3)&=-\frac{7 \zeta (3)}{\pi ^2}+\frac{310 \zeta (5)}{\pi ^4}-\frac{1905 \zeta (7)}{\pi ^6},\\
J(5,5)&=\frac{31 \zeta (5)}{\pi ^4}-\frac{3175 \zeta (7)}{\pi ^6}+\frac{35770 \zeta (9)}{\pi ^8},\\
J(6,2)&=\frac{322 \zeta (3)}{15 \pi ^2}-\frac{248 \zeta (5)}{\pi^4}+\frac{762 \zeta (7)}{\pi ^6},\\
J(6,4)&=-\frac{2852 \zeta (5)}{15 \pi ^4}+\frac{5080 \zeta (7)}{\pi ^6}-\frac{28616 \zeta (9)}{\pi ^8},\\
J(6,6)&=\frac{5842 \zeta (7)}{5 \pi ^6}-\frac{57232 \zeta (9)}{\pi^8}+\frac{515844 \zeta (11)}{\pi ^{10}},\\
J(7,7)&=-\frac{127 \zeta (7)}{\pi ^6}+\frac{1402184 \zeta (9)}{45 \pi^8}
	-\frac{1003030 \zeta (11)}{\pi ^{10}}+\frac{7568484 \zeta (13)}{\pi ^{12}}.}

Among these formulae, the first one was proposed as a monthly problem
by Holland~\cito{holl}, which has been the primary motivation
for the authors to carry on the present research.


\end{document}